\documentclass[11pt,twoside]{amsart} \textwidth 12 cm
\numberwithin{equation}{section}

\begin{document}

\title{A note on the possibility  of a  motion  without crossing a given subset}
\author{R. Mirzaie}

\begin{abstract}   We study  the frcatl dimension
of a given subset $X$ of $R^{n}$ such that a motion is possible without crossing $X$.
\end{abstract}
\thanks{{\scriptsize
\flushleft 1. Department of Mathematics, Faculty of Science, Imam Khomeini International University (IKIU), Qazvin, Iran.
 Email: r.mirzaei@sci.ikiu.ac.ir
{\scriptsize
}}}
\maketitle

\pagestyle{myheadings}

\markboth{\rightline {\scriptsize  R. Mirzaie}}
         {\leftline{When there is short  way to a submanifold }}
\section{Introduction}
Let $X$ and $Y$ be subsets of $R^{n}$. Suppose that the points of $X$ are spread all over $R^{n}$ and consider $Y$ as an object
which we intend to move it in $R^{n}$  without crossing $X$. Study of the possibility
of a motion under given conditions is a general problem   which
have many applications in physics and other areas of science.
Possibility of the  mentioned motion of $Y$  is dependent to the
size of $X$ and $Y$. For a mathematical study, we
 must mention what we mean from  "motion" and
"size". In the present note,  we define a "motion" in $R^{n}$.
Then, we use the fractal dimension of  objects as their size to
study of the motion  possibility. Also, we show that our main result can easily be generalized to the Riemannian manifolds\\\\
The following definition is motivated from the definition of
homotopy through
isometries in geometry.\\\\
 {\bf Definition 1.1.} Let $Y \subset
R^{n}$. A motion of $Y$ in $R^{n}$ is a continuous map $\mu:Y
\times [0,1] \to R^{n}$ such that $\mu(y,0)=y$, $y \in Y$,  and
for all $t \in [0,1]$, the following map is an isometry:
\[ Y \to \mu(Y,t), \ \ \ \ y \to \mu(y,t)\]
 $Y_{1}=\mu(Y,1)$ is called the target of the motion, and for a fixed $y$, the following map is called
 the path of $y$
 \[ \mu_{y}:[0,1] \to R^{n}, \ \ \ \mu_{y}(t)= \mu(y,t).\]
 There are many candidates to use as the size of a subset of
 $R^{n}$. Two important candidates are \\
 (1) Measure of the set by suitable definition of a measure on
 $R^{n}$.\\
 (2) Dimension of the set ( using several definitions of
 dimension, like hausdroff dimension, box dimension, and so on).\\
 In this note, we will use the box dimension, which we remind its
 definition.
Let $A$ be a subset of a metric space $(M, d)$. We denote by
dim$A$ the topological dimension of $A$. If $A$  is bounded then
the upper box dimension of $A$ is defined by
\[ \overline{dim}_{B}(A)= limsup_{\delta \to
0}\frac{m_{\delta}(A)}{-log \delta}\] where, $m_{\delta}A$ is the
maximum number of disjoint balls of radius $\delta$, with centers
contained in $A$. The lower box dimension $\underline{dim}_{B}(A)$
is defined in similar way. If the upper and the lower box dimensions  of
$A$ are equal, box dimension of $A$ is defined (see [1] for details, and other definitions of dimension). In what follows we suppose that the box dimension
exists. If not, we can use upper or lower box dimensions every where
.\\
If $A$ is not bounded  then the box dimension of $A$ is defined by
\[ dim_{B}A = max \{ dim_{B}(A \cap U): \ \  U  \ \ is \ \ open \ \  bonded \ \ set \}
 .\]
{\bf Problem 1.2.} { Let $ Y \subset R^{n}$ and $s \in R^{n}, s
\notin  Y$. Suppose that $s$ is a light source and the whole space
is covered by a dust which we denote it by $X$. How much big can
be
 dim$_{B}X$ such that $Y$ is not covered by  the shadow of $X$.}\\\\
{\bf Problem 1.3.}  Consider $X$ as Problem 1.2, and let $M$ be a
submanifold of $R^{n}$.  How much big can be dim$_{B}X$ such that
$M$ can move in
$R^{n}$, without crossing $X$.\\\\
Answers of the problems (1.2) and (1.3) are equivalent to the
following theorems (2.1) and (2.2) related to the motion
possibility problem.
\section {ResultsAAA}

 {\bf Theorem 2.1.} {\it Let $Y, X$ be a closed subsets of $R^{n}$, $n \geq
 2$, such that dim$_{B} X <$dim$_{B}Y$.
   Then for each point $s \in R^{n}-X$ there exists a motion of $\{s\}$ such that its target is included in $Y$
   and all paths  are straight lines.
   }\\\\
   {\bf Proof:}
 Since $X$ is closed, there exists $r>0$
such that
   $B(s,r) \cup X=\emptyset$. Consider the following set, which is
   called the cone with vertex $s$ over $Y$:
   \[ C(s,Y)=\{ ts+(1-t)y: \ \ 0\leq t \leq 1, \ \ \ , y \in Y\}\]
   Consider the following map
   \[ f: C(s,Y)-B(r,s) \to Y, \ \ \  f(ts+(1-t)y)=y\]
   We show that $f$ is Lipschitz. \\
   Let $a_{1}=t_{1}s +(1-t_{1})y_{1}$ and
   $a_{2}=t_{2}s+(1-t_{2})y_{2}$.  Consider the triangle
   $sy_{1}y_{2}$ and let $s_{1}$ and $s_{2}$ be the intersection
   points of the line segments $sy_{1}$ and $sy_{2}$ with the
   boundary of $B(s,r)$. By Thalis equality we have
   \[
   \frac{d(y_{1},y_{2})}{d(a_{1},a_{2})}=\frac{d(y_{1},s)}{d(a_{1},s)} \ \ \ \ (1)\]
   Since $a_{1} \notin B(s,r)$, then $d(a_{1},s)>d(s_{1},s)=r$.
   Thus by (1)
   \[\frac{d(y_{1},y_{2})}{d(a_{1},a_{2})}<\frac{d(y_{1},s)}{r} \ \ \ \ (2)\]
   Without lose of generality we can assume that $Y$ is compact.
   Then there exists a maximum for  the set  of numbers $\{
   d(y_{1},s): y_{1} \in Y\}$, which we denote it by $R$. Then we
   get from (2)
   \[\frac{d(y_{1},y_{2})}{d(a_{1},a_{2})}<\frac{R}{r}\]
   Then we have
   \[\frac{d(y_{1},y_{2})}{d(a_{1},a_{2})}<\frac{R}{r} \Rightarrow
   \frac{d(f(a_{1}),f(a_{2}))}{d(a_{1},a_{2})}<\frac{R}{r}\]
   So, $f$ is Lipschitz. Now we have
   \[ dim_{B}(f(X \cap C(Y,s)))\leq dim_{B}(X) <dim_{B} Y\]
   Therefore, $Y-Im(f)$ in not empty and there exists a point $y
   \in Y$ such that $y$ does not belong to the image of $f$. So
   the line segment $ts+(1-t)y$, $0\leq t\leq 1$, does not
   intersect $X$.\\\\
{\bf Theorem 2.2.} {\it \\(1) If $M$ is a submanifold of $R^{n}$
and
   dim$X<n-dimM-1$ and $M \cap X=\emptyset$, then  there exists a
   a motion of $M$ such that all paths are disjoint from $X$.\\
   (2) Given $x_{0} \in M$ and $y_{0} \in R^{n}-X$, we can choose the motion in (1) in such a way that the target of $x_{0}$ be equal to
   $y_{0}$.
    }\\\\
{\bf Proof:} Step 1. If $K \subset R^{n}$ and dim$_{B} K <n-1$,
then there exists a path $\alpha:[0,1] \to R^{n}$ such that
$\alpha(0,1]
\cap K =\emptyset$.\\
 Proof:  Consider a point $b_{1} \in R^{n}-K$, with $|b_{1}|>
1$. If we put in Theorem 1, $X=S^{n-1}-K$, $M=K$, then we get that
there exists a motion of $b_{1}$ such that its destination is
included in $S^{n-1}(1)-K$ and the path of the motion is a line
segment. Thus, there is a point $b_{2} \in S^{n-1}(1) $ such that
the line segment $b_{1}b_{2}$ does not intersect $K$. In a
similar way, there exists a point $b_{3} \in
S^{n-1}(\frac{1}{2})$ such that the line segment $b_{2}b_{3}$
does not intersect $K$. If we proceed in this we we obtain
sequence of line segments $b_{1}b_{2},
b_{2}b_{3},...,b_{m}b_{m+1},...$, such that $|b_{m}| \to 0$. By a
suitable choice of parameter, we can define a path $\alpha:[0,1]
\to R^{n}$, such that $\alpha(0)=\circ$ and the image of $\alpha$
is union of the origin of $R^{n}$ and the
mentioned line segments.\\
Step 2. Consider the map $f: M \times R^{n} \to R^{n}$
 defined by $f(x,y)=y-x$. Put $K=f(M \times X)$.
 Consider the following metric on $M \times R^{n}$
 \[ d((x_{1},y_{1}),(x_{2},y_{2}))=|x_{1}-x_{2}|+|y_{1}-y_{2}|, \ \ (x_{i},y_{i}) \in X \times R^{n}\]
 Clearly $f$ is a Lischitz map, so
 \[ dim_{B}(K)=dim_{B}f(M\times X) \leq dim_{B}(M)+dim_{B}(X)<n-1\]
 So, by step 1, there exists a path $\alpha:[0,1] \to R^{n}$ such that
 $\alpha(0,1] \cap K=\emptyset$. Put
 \[ F: M \times [0,1] \to R^{n},  F(x,t)=x+\alpha(t)\]
 $F$ is a homotopic motion which does not intersect $X$.
 Because, if for some $x \in M$ and $t \in (0,1]$,  $F(x,t) \in
 Y$, then
 \[ (F(x,t)-x) \in K \Rightarrow \alpha(t) \in K\]
 which is contradiction. This proves (1).\\
 For the proof of (2), consider the point $b_{1}$ in the proof of
 (1), which is the end point of $\alpha$. $b_{1}$ is arbitrary
 point in $R^{n}-Y$. So we can choice it as $b_{1}=y_{0}-x_{0}$.
 Then we have
 \[ F(x_{0},1)=x_{0}+\alpha(1)=y_{0}.\]\\\\
 {\bf Remark 2.3.} In Theorem 2, we can replace $R^{n}$ by any open subset of $R^{n}$ containing $X$.\\\\
    {\bf Proof:} In the the proof of Theorem 2, Step 1, we have $|\alpha(t)|\leq |b_{1}|$, for all $t \in [0,1]$. Since dim$_{B}(K) <n-1$, then $b_{1}$ can be chosen arbitrary close to the origin. So, $|\alpha(t)|$ can
 be arbitrary small. Then $|F(x,t)-x|$ can be arbitrary small.\\\\
 {\bf Remark 2.4 (Optimality).}
 It is easy to give examples to show that upper bound estimates for
 dim$_{B}X$ in the theorems 1 and 2 are optimal. \\
1.  If $X=\{(x_{1}, ..., x_{n-1},0), x_{i} \in R \}$ and
$Y=\{(x_{1}, ..., x_{n-1},1), x_{i} \in R\}$ and
$s=(1,1,1...,1,2)$, there
 is no motion of $\{b\}$ with targets on $X$ without crossing
 $Y$.\\
 2. Let $Q$ denote the rational numbers. Put $n=3$, $M=\{(0,0,z): z \in R\}$,  $X=Q^{3}-M$.  Since $X$
 is countable then dim$_{B}(X)=0$, so dim$_{B}(X)=n-dim_{B}(M)-1$.
 Any nontrivial motion of $M$ cuts $X$.\\\\
  Mentioned problems can be be considered in more general case
    where $R^{n}$ is replaced by a Riemannian manifold $N$. A motion in $N$ is defined in a similar way as definition 1.1, replacing
    $R^{n}$ by $N$.  We prove the
    following theorem which is generalization of Theorem 2.2.\\\\
    {\bf Theorem 2.5.} {\it Let $N$ be a Riemannian manifold of dimension $n$, $M$ be a submanifold of $N$ and
   $Y$ be a subset of $M$ such that dim$_{B}Y<n-dimM-1$ and $M \cap Y=\emptyset$, then  there exists
   a motion of $M$ such that all paths are disjoint from $Y$.\\
   (2) Given $x_{0} \in M$ and $d \in M-Y$, we can choose the motion in (1) in
   such a way that the target of $x_{0}$ be equal to
   $d$.  }\\\\
 {\bf Proof:}
By Nash's embedding theorem, for a sufficiently big integer
 number $m$, $N$ can be considered as a submanifold of
 $R^{m}$.

 For each point $p \in N$ there is a neibourhood $W$ of $p$ in $N$ and a positive integer $\epsilon=\epsilon(p)$
 such that the set
  \[\hat{W}(p, \epsilon(p))=\{ x+tv: x \in W, v \in (T_{p}N)^{\perp}, |v|=1, t \in (- \epsilon, \epsilon)\}\]
   is an open set in $R^{m}$. Put
 \[ \hat{N}=\bigcup_{p \in N} \hat{W}(p, \epsilon(p)), \ \ \hat{Y}=\bigcup_{p \in Y} \hat{W}(p, \epsilon(p))\]
 Clearly dim$_{B}(\hat{Y}) \leq$ dim$_{B}Y +(m-n)$. So,  dim$_{B}(\hat{Y})  <(n-dimM-1)+(m-n)=m-dim M-1$.
 Then by Theorem 2.2 and Remark 2.3, there is a motion $\hat{\mu}: M \times [0,1] \to R^{m}$ which does not intersect $\hat{Y}$,
 that is $\hat{\mu}(M \times (0,1]) \bigcap \hat{Y} =\emptyset$. Now, put
 \[ \pi :\hat{N} \to N,  \pi(x+tv)=x\]
 and
 \[ \mu:M\times [0,1] \to N, \ \  \mu(x,t)=\pi \circ \hat{\mu}(x,t)\]
 $\mu$ is the desired motion.

\end{document}